# The Probabilities of Large Deviations for the Chi-square and Log-likelihood Ratio Statistics

Sherzod Mirakhmedov

Institute of Mathematics. National University of Uzbekistan
100125, Tashkent, Durmon yuli st., 29
e-mail: shmirakhmedov@yahoo.com

**Abstract**. A new large deviation results for the Pearson's chi-square and Log-likelihood ratio statistics are obtained. Here attention is focused on the case when the number of groups increases to infinity and the probabilities of groups' decreases to zero, as the sample size tends to infinity.

**1. Introduction**. Let $\eta = (\eta_1,...,\eta_N)$ be a random vector having a multinomial distribution $M(n, p_1,...,p_N)$, where $p_1 +...+ p_N = 1$, $p_m > 0$, $m = 1,...,N$ and $\eta_1 +...+ \eta_N = n$. Many tests of goodness-of-fit can be reduced to testing a hypothesis about the parameters $p_1,...,p_N$. Most notably test statistics for this kind of testing problem are Pearson's chi-square ($\chi^2$) statistic

$$\chi_N^2 = \sum_{m=1}^{N} \frac{(\eta_m - np_m)^2}{np_m} ,$$

and the log-likelihood ratio (LR) statistic:

$$\Lambda_N = 2\sum_{m=1}^{N} \eta_m \ln \frac{\eta_m}{np_m} .$$

In the present paper we study $P\{T_N \geq x_N\}$, where $T_N$ is standardized version of the $\chi_N^2$ or $\Lambda_N$ statistics and $x_N \to \infty$, in the situation when $\max p_m \to 0$, $N \to \infty$ as $n \to \infty$.

Similar problems have been considered in the literature. We refer to Kallenberg (1985), in particularly his results for the Cramer's zone, viz. $x_N = o(\sqrt{N})$, proposes that $N \min p_m \geq c > 0$ and additionally $N = o(\sqrt{n})$, $\sqrt{n} N^{-3/2} x_N \to \infty$ for the $\chi_N^2$ case. We refer also to Sirajdinov et al (1989) and Ivchenko and Mirakhmedov (1995). In these papers the large deviation problem in Cramér's zone has been studied for the class of statistics $h_1(\eta_1) +...+ h_N(\eta_N)$, where $h_m(u)$ is a Borel function defined on the non-negative axis and satisfying the Cramer condition: $\max_{1 \leq k \leq N} E \exp\{H|h_k(\xi)|\} < \infty$, $\exists H > 0$, here $\xi \sim Poi(\lambda)$, the Poisson distribution with expectation $\lambda = \lim_{n \to \infty}(n/N) \in (0,\infty)$. The Cramér's condition is fulfilled for the $\Lambda_N$, but for the $\chi_N^2$ does not

**MSC (2000):** 60F10, 62G10, 62G20

**Key words and phrases**. Chi-square statistic, Log-likelihood ration statistic, multinomial distribution, large deviations, Poisson distribution.

---

The assertions of the present paper corresponds to the case $x_N = o(N^{1/6})$, nevertheless they covers the situations which does not covered by aforementioned papers. Our results, together with corresponding results of Kallenberg (1985), allowed to explore the intermediate asymptotic relative efficiency due to Inglot (1999) of the $\chi^2$ and LR goodness of fit tests; the results are object of separate paper Mirakhmedov (2016). Also our proofs are based on the technique different from those aforementioned authors. Core part of our technique is that we reduce the cumulants of statistics of type $R_N$ to the cumulants of a sum of independent random variables (r.v.), next we use a modified variant of Lemma 2.3 of Saulis and Statulevicius (1991).

The rest of the paper is organized as follows. The main results are presented in Section 2; in Section 3 the proofs are given; for the reader's convenience, the auxiliary Assertions are collected in Appendix. In what follows $C_k$ is an universal positive constant; all asymptotic statements are considered as $n \to \infty$; $\varsigma \sim F$ stands for " r.v $\varsigma$ has the distribution $F$ ".

## 2. Results

Let $p_{\max} = \max_{1 \leq k \leq N} p_k$, $p_{\min} = \min_{1 \leq k \leq N} p_k$, $\lambda_n = n/N$ and $\Phi(u)$ be the standard normal distribution function.

**2.1. Chi-square statistic**. Set $\nabla_n = \max\left(1, (np_{\min})^{-1}\right)$,

$$\tilde{\sigma}_N^2 = \sum_{m=1}^N \frac{1}{np_m} + 2N, \quad \sigma_N^2 = \frac{1}{\lambda_n}\sum_{m=1}^N \left(\frac{1}{Np_m} - 1\right) + 2N = \tilde{\sigma}_N^2 - N\lambda_n^{-1}. \tag{2.1}$$

**Theorem 2.1.** Let $p_{\max} = o(1)$. Then for all $x_N$ such that $0 \leq x_N = o\left(\min\left(\left(\sigma_N^3/\tilde{\sigma}_N^2 \nabla_n\right)^{1/3}, n^{1/6}, p_{\max}^{-1/4}\right)\right)$ it holds

$$P\{\chi_N^2 > x_N \sigma_N + N\} = (1 - \Phi(x_N))(1 + o(1)) \tag{2.2}$$

and

$$P\{\chi_N^2 < -x_N \sigma_N + N\} = \Phi(-x_N)(1 + o(1)). \tag{2.3}$$

**Corollary 2.1.** Let $\lambda_n \geq C_0 > 1/2$ and

$$C_1 \leq Np_{\min} \leq Np_{\max} \leq C_2 N^{1/3}, \tag{2.4}$$

then for all $x_N$ such that $0 \leq x_N = o(N^{1/6})$ the relations (2.2) and (2.3) are valid.

**Corollary 2.2.** Let $p_m = N^{-1}\left(1 + \delta(n)d_{m,n}\right), m = 1,...,N$, where $\delta(n) \to 0$ and

$$\sum_{m=1}^N d_{m,n} = 0, \quad \frac{1}{N}\sum_{m=1}^N d_{m,n}^2 = d^2 < \infty.$$

Then for arbitrary $\lambda_n$ and $x_N$ such that $0 \leq x_N = o\left(N^{1/6} \min(1, \lambda_n^{2/3})\right)$ the relations (2.2) and (2.3) are valid.

**Remarks 2.1**. For $x_N \to \infty$ the relation (2.2) can be written in the form $\log P\{\chi_N^2 > x_N \sigma_N + N\} = -x_N^2/2(1+o(1))$. From Kallenberg (1985, Eq. (2.17)) such kind relation follows if $\log N = o(x_N^2)$, $x_N = o(N^{1/6})$ and $N = o(n^{3/8})$. Last condition excludes the case recommended by Mann and Wald (1942) who obtained the relation $N \sim cn^{2/5}$ concerning the optimal choice of the number of groups in chi-square goodness of fit test.

**2.2. Likelihood ratio statistic $\Lambda_N$**. Here we restrict our attention to the case $\lambda_n \to \infty$; the corresponding results for the case $\lambda_n \to \lambda \in (0, \infty)$ follows from Sirajdinov et al (1987).

**Theorem 2.2.** Let $\lambda_n \to \infty$ and

$$C_3 \leq Np_{\min} \leq Np_{\max} \leq C_4 \lambda_n^3, \tag{2.5}$$

Then for all $x_N$ such that $0 \leq x_N = o\left(\min(N^{1/6}, p_{\max}^{-1/4})\right)$ it holds

$$P\{\Lambda_N > x_N \sqrt{2N} + N\} = (1 - \Phi(x_N))(1+o(1)), \tag{2.6}$$

and

$$P\{\Lambda_N < -x_N \sqrt{2N} + N\} = \Phi(-x_N)(1+o(1)). \tag{2.7}$$

**Corollary 2.3.** Let $\lambda_n \to \infty$ and $p_m = N^{-1}(1+o(1))$, $m=1,...,N$, then for all $x_N$ such that $0 \leq x_N = o(N^{1/6})$ the relations (2.6) and (2.7) are valid.

**Remarks 2.2.** Let $\lambda_n \to \infty$, $Np_{\max} = O(\min(N^{1/3}, \lambda_n^3))$, $x_N \to \infty$ and $x_N = o(N^{1/6})$. Then from (2.6) it follows that $\log P\{\Lambda_N > x_N \sqrt{2N} + N\} = -x_N^2/2(1+o(1))$. From Kallenberg (1985, Eq. (2.13)) such relation follows under additionally conditions $N = o(n^{3/7})$ and $\log N = o(x_n^2)$.

### 3. Proofs

In what follows $\xi_1, ..., \xi_N$ are independent r.v.s with $\xi_m \sim Poi(np_m)$, $\xi \sim Poi(\lambda_n)$ and $\mathcal{C}_k(\zeta)$ is a cumulant of the $k$ th order of the r.v. $\zeta$,

$$R_N = \sum_{m=1}^{N} h_m(\eta_m) \quad \text{and} \quad T_N = \sum_{m=1}^{N} h_m(\xi_m), \tag{3.1}$$

where $(\eta_1, ..., \eta_N)$ is a multinomial r.vec. as in introduction, $h_m(x)$, $m=1,...,N$ are Borel functions defined on the non-negative axis.

**Main Lemma.** Let

$$E|h_m(\xi_m)|^s = O\left(Eh_m^s(\xi_m)\right), \tag{3.2}$$

$$E|h_m(\xi_m)|^s (\xi_m - np_m)^2 = O\left(s^a (np_m)^b Eh_m^s(\xi_m)\right), \tag{3.3}$$

for all $m = 1, 2, ..., N$, some non-negative real $a$ and $b$, and each integer $s: 3 \leq s \leq k$, where an integer $k$ such that

$$\sum_{i=1}^{k} p_{j_i} = o(1), \tag{3.4}$$

for any $k$-subset $(j_1, ..., j_k)$ of the set $(1, 2, ..., N)$. Then for all $k$ satisfying (3.4) and

$$3 \leq k = o(K_n(a,b)), \tag{3.5}$$

where $K_n(a,b) = \left(n^{1-b} p_{\max}^{-b}\right)^{1/(\bar{a}+1)}$, $\bar{a} = \max(1, a)$, it holds

$$\mathcal{C}_k(R_N) = \mathcal{C}_k(T_N)(1 + o(1)). \tag{3.6}$$

**Proof.** We start from the following

**Lemma 3.1.** For any non-negative $k$ it holds

$$ER_N^k = \upsilon_n \int_{-\pi\sqrt{n}}^{\pi\sqrt{n}} ET_N^k \exp\left\{i\tau \frac{S_N}{\sqrt{n}}\right\} d\tau,$$

where

$$\upsilon_n \stackrel{def}{=} \left(2\pi\sqrt{n} P\{S_N = 0\}\right)^{-1} = \frac{n! e^n}{2\pi n^n \sqrt{n}} = \frac{1}{\sqrt{2\pi}}\left(1 + o\left(\frac{1}{n}\right)\right). \tag{3.7}$$

**Proof.** It is well known that $\mathcal{L}((\eta_1, ..., \eta_N)) = \mathcal{L}((\xi_1, ..., \xi_N)/S_N = 0)$, where $\mathcal{L}(X)$ stands for the distribution of a random vector $X$, $S_N = \sum_{m=1}^{N}(\xi_m - np_m)$. Hence $ER_N^k = E(T_N^k | S_N = 0)$. On the other hand $E(T_N^k e^{i\tau S_N}) = E\{e^{i\tau S_N} E(T_N^k | S_N)\}$. Now Lemma 3.1 follows by Fourier inversion. Notice that Eq. (3.7) follows due to fact that $S_N + n \sim Poi(n)$ and Stirling's formula.

Let $k$ be as in (3.4) and (3.5). Use Lemma 3.1 to write

$$ER_N^k = \mu_n \sum_{l=1}^{k} \sum_{l,k}' \sum_{l}'' \int_{-\pi\sqrt{n}}^{\pi\sqrt{n}} E\left(\left(h_{j_1}(\xi_{j_1})\right)^{k_1} \cdot ... \cdot \left(h_{j_l}(\xi_{j_l})\right)^{k_l} \exp\left\{i\tau \frac{S_N}{\sqrt{n}}\right\}\right) d\tau, \tag{3.8}$$

where $\sum_{l,k}'$ is the summation over all $l$-tuples $(k_1, ..., k_l)$ with positive integer components such that $k_1 + ... + k_l = k$; $\sum_{l}''$ is the summation over all $l$-tuples $(j_1, ..., j_l)$ such that $j_i \neq j_r$ for $i \neq r$ and $j_m = 1, 2, ..., N$; $m = 1, 2, ..., l$.

Set $S_{l,N} = \sum_{i=1}^{l}(\xi_{j_i} - np_{j_i})$, $d_{l,N} = \sum_{i=1}^{l} p_{j_i}$, $l = 1, ..., k$ and write

$$\int_{-\pi\sqrt{n}}^{\pi\sqrt{n}} E\left(\left(h_{j_1}(\xi_{j_1})\right)^{k_1} \cdot ... \cdot \left(h_{j_l}(\xi_{j_l})\right)^{k_l} \exp\left\{i\tau \frac{S_N}{\sqrt{n}}\right\}\right) d\tau$$

$$= \int_{-\pi\sqrt{n}}^{\pi\sqrt{n}} E\left[\left(h_{j_1}(\xi_{j_1})\right)^{k_1} \cdot \ldots \cdot \left(h_{j_l}(\xi_{j_l})\right)^{k_l} \exp\left\{i\tau \frac{S_{l,N}}{\sqrt{n}}\right\}\right] E \exp\left\{i\tau \frac{S_N - S_{l,N}}{\sqrt{n}}\right\} d\tau$$

$$= \int_{-\pi\sqrt{n}}^{\pi\sqrt{n}} E\left[\left(h_{j_1}(\xi_{j_1})\right)^{k_1} \cdot \ldots \cdot \left(h_{j_l}(\xi_{j_l})\right)^{k_l} \exp\left\{i\tau \frac{S_{l,N}}{\sqrt{n}}\right\}\right] \left(E \exp\left\{i\tau \frac{S_N - S_{l,N}}{\sqrt{n}}\right\} - \exp\left\{-\frac{\tau^2}{2}(1-d_{l,N})\right\}\right) d\tau$$

$$+ \int_{-\pi\sqrt{n}}^{\pi\sqrt{n}} \exp\left\{-\frac{\tau^2}{2}(1-d_{l,N})\right\} E\left[\left(h_{j_1}(\xi_{j_1})\right)^{k_1} \cdot \ldots \cdot \left(h_{j_l}(\xi_{j_l})\right)^{k_l} \left(\exp\left\{i\tau \frac{S_{l,N}}{\sqrt{n}}\right\} - 1\right)\right] d\tau$$

$$+ E\left[\left(h_{j_1}(\xi_{j_1})\right)^{k_1} \cdot \ldots \cdot \left(h_{j_l}(\xi_{j_l})\right)^{k_l}\right] \int_{-\pi\sqrt{n}}^{\pi\sqrt{n}} \exp\left\{-\frac{\tau^2}{2}(1-d_{l,N})\right\} d\tau \stackrel{def}{=} J_1 + J_2 + J_3. \tag{3.9}$$

We have
$$E\exp\left\{\frac{i\tau(\xi_m - np_m)}{\sqrt{n}}\right\} = \exp\left\{np_m\left(e^{it/\sqrt{n}} - 1 - \frac{it}{\sqrt{n}}\right)\right\} = \exp\left\{-\frac{\tau^2}{2}p_m + \frac{\theta\tau^3}{6\sqrt{n}}p_m\right\}, \tag{3.10}$$

and
$$\left|E\exp\left\{\frac{i\tau(\xi_m - np_m)}{\sqrt{n}}\right\}\right| = \exp\left\{-2np_m \sin^2 \frac{\tau}{2\sqrt{n}}\right\} \leq \exp\left\{-\frac{2p_m}{\pi^2}\tau^2\right\}, \tag{3.11}$$

because $\sin^2 u/2 \geq u^2/\pi^2$, $|u| \leq \pi$.

In order to get an upper bound for the $|J_1|$ we first write the integral $J_1$ as the sum of two integrals, $J_1'$ and $J_1''$ say, over intervals $|\tau| \leq \pi\sqrt{n}/2$ and $\pi\sqrt{n}/2 \leq |\tau| \leq \pi\sqrt{n}$ respectively. Next, noting that $S_N - S_{l,N}$ is the sum of independent r.v.s $\xi_m - np_m$ over the set of indexes $(1,...,N)$ except $(j_1,...,j_l)$ we use to $E\exp\{i\tau(S_N - S_{l,N})/\sqrt{n}\}$ in the $J_1'$ and $J_1''$ the equation (3.10) and inequality (3.11) respectively. Then quite clear algebra gives

$$|J_1| \leq E|h_{j_1}(\xi_{j_1})|^{k_1} \cdot \ldots \cdot E|h_{j_l}(\xi_{j_l})|^{k_l} \left\{\int_{-\pi\sqrt{n}/2}^{\pi\sqrt{n}/2} \exp\left\{-\frac{\tau^2}{2}(1-d_{l,N})\right\} \left(\exp\left\{\frac{|\tau|^3}{6\sqrt{n}}(1-d_{l,N})\right\} - 1\right) d\tau\right.$$

$$+ \int_{\frac{\pi\sqrt{n}}{2} \leq |\tau| \leq \pi\sqrt{n}} \left(\left|E\exp\left\{i\tau \frac{S_N - S_{l,N}}{\sqrt{n}}\right\}\right| + \exp\left\{-\frac{\tau^2}{2}(1-d_{l,N})\right\}\right) d\tau\right\}$$

$$\leq E|h_{j_1}(\xi_{j_1})|^{k_1} \cdot \ldots \cdot E|h_{j_l}(\xi_{j_l})|^{k_l} \left[\int_{-\pi\sqrt{n}/2}^{\pi\sqrt{n}/2} \exp\left\{-\frac{6-\pi}{12}(1-d_{l,N})\tau^2\right\} \frac{|\tau|^3}{6\sqrt{n}}(1-d_{l,N}) d\tau\right.$$

$$\left. + C_5 \exp\left\{-\frac{n}{4}(1-d_{l,N})\right\}\right].$$

Hence
$$J_1 = O\left(\frac{1}{\sqrt{n}}\right) E\left(\left(h_{j_1}(\xi_{j_1})\right)^{k_1} \cdot \ldots \cdot \left(h_{j_l}(\xi_{j_l})\right)^{k_l}\right). \tag{3.12}$$

since (3.2) and (3.4). Next apply inequality $|e^{it} - 1 - it| \leq t^2/2$ to get: for all $k$ satisfying (3.5)

$$J_2 = \frac{C_6}{n} E\left( |h_{j_1}(\xi_{j_1})|^{k_1} \cdot \ldots \cdot |h_{j_l}(\xi_{j_l})|^{k_l} \left(\sum_{i=1}^{l}(\xi_{j_i} - np_{j_i})\right)^2 \right)$$

$$\leq \frac{C_7}{n} l \sum_{i=1}^{l} E|h_{j_i}(\xi_{j_i})|^{k_i} (\xi_{j_i} - np_{j_i})^2 \prod_{m \neq i, m=1}^{l} |h_{j_m}(\xi_{j_m})|^{k_m}$$

$$\leq C_8 E\left( (h_{j_1}(\xi_{j_1}))^{k_1} \cdot \ldots \cdot (h_{j_l}(\xi_{j_l}))^{k_l} \right) \frac{k}{n} \sum_{i=1}^{l} k_i^a (np_{j_i})^b$$

$$\leq C_9 E\left( (h_{j_1}(\xi_{j_1}))^{k_1} \cdot \ldots \cdot (h_{j_l}(\xi_{j_l}))^{k_l} \right) \frac{(np_{\max})^b}{n} k^{1+\max(1,a)}$$

$$= o(1) E\left( (h_{j_1}(\xi_{j_1}))^{k_1} \cdot \ldots \cdot (h_{j_l}(\xi_{j_l}))^{k_l} \right), \quad (3.13)$$

because (3.3) and the fact that $k_1^a + \ldots + k_l^a = k^{\max(1,a)}$. By a simple algebra we obtain

$$J_3 = \sqrt{2\pi} E\left( (h_{j_1}(\xi_{j_1}))^{k_1} \cdot \ldots \cdot (h_{j_l}(\xi_{j_l}))^{k_l} \right)(1 + o(1)). \quad (3.14)$$

since the condition (3.4). Now apply (3.12), (3.13) and (3.14) in the (3.9) to get

$$\int_{-\infty}^{\infty} E\left( (h_{j_1}(\xi_{j_1}))^{k_1} \cdot \ldots \cdot (h_{j_l}(\xi_{j_l}))^{k_l} \exp\left\{ i\tau \frac{S_N}{\sqrt{n}} \right\} \right) d\tau$$

$$= \sqrt{2\pi}(1 + o(1)) E\left( (h_{j_1}(\xi_{j_1}))^{k_1} \cdot \ldots \cdot (h_{j_l}(\xi_{j_l}))^{k_l} \right). \quad (3.15)$$

Remark that in (3.13) and (3.14), and hence in (3.15), the $o(1)$ is uniform in all $l$-tuples $(k_{j_1}, \ldots, k_{j_l})$ and $k$ satisfying the condition (3.5). The relations (3.7), (3.8) and (3.15) imply

$$ER_N^k = ET_N^k (1 + o(1)), \quad (3.16)$$

for every integer $k$ satisfying the condition (3.5). Main Lemma follows from (3.16) and Assertion 3.

**Remark 3.1**. The condition (3.3) is used to estimate of the integral $J_2$ only. Let the condition (3.3) is fulfilled with $a_1, b_1$ for all $m \in \mathbb{N} \subset (1, 2, \ldots, N)$ and with $a_2, b_2$ for all $m \in (1, 2, \ldots, N) \setminus \mathbb{N}$, then Eq. (3.6) is valid for every integer $k$ such that $3 \leq k = o\left(\min(K_n(a_1, b_1), K_n(a_2, b_2))\right)$ instead of (3.5).

Let $\lfloor x \rfloor$ stand for the integer part of $x$, $\xi \sim Poi(\lambda)$ and $\mu_v(\lambda) = E(\xi - \lambda)^v$.

**Lemma 3.2.** For any integer $v \geq 2$ one has

$$\mu_v(\lambda) = v! \sum_{l=1}^{\lfloor v/2 \rfloor} c_{l,v} \lambda^l, \quad (3.17)$$

where for all of $l = 1, 2, \ldots, \lfloor v/2 \rfloor$,

$$0 < c_{l,v} < 1/l!, \quad (3.18)$$

also

$$(v+1)c_{l,v+1} = lc_{l,v} + c_{l-1,v-1}, \qquad (3.19)$$

where $l = 1, 2, ..., \lfloor (v+1)/2 \rfloor$ if $v$ is an even, and $l = 1, 2, ..., \lfloor (v+1)/2 \rfloor - 1$,

$(v+1)c_{\lfloor (v+1)/2 \rfloor, v+1} = c_{\lfloor (v-1)/2 \rfloor, v-1}$ if $v$ is an odd; here we put $c_{0,v-1} = 0$.

**Proof.** Apply the Bruno's formula to $Ee^{i\tau(\xi-\lambda)} = \exp\{\lambda(e^{i\tau} - 1 - i\tau)\}$ to get

$$\mu_v(\lambda) = v! \sum \lambda^{k_2+...+k_v} \prod_{m=2}^{v} \frac{1}{k_m!(m!)^{k_m}} = v! \sum_{l=1}^{\lfloor v/2 \rfloor} c_{l,v} \lambda^l,$$

where $\sum$ is summation over all non-negative $k_2, ..., k_v$ such that $2k_2 + ... + vk_v = v$ and

$l = k_2 + ... + k_v$, $c_{l,v} = \sum \prod_{m=2}^{v} \frac{1}{k_m!(m!)^{k_m}}$, and hence $0 < c_{l,v} < 1/l!$. On the other hand Kenney and Keeping (1953) give

$$\mu_{v+1}(\lambda) = v\lambda \mu_{v-1}(\lambda) + \lambda \frac{d}{d\lambda} \mu_v(\lambda). \qquad (3.20)$$

Just to keep notation simple we put $\kappa = \lfloor v/2 \rfloor$ and $\hat{c}_{i,v} = v!c_{i,v}$. Note that $\hat{c}_{1,v} = 1$ for any integer $v \geq 0$. Let $v$ is even, then $\lfloor (v+1)/2 \rfloor = \kappa$ and $\lfloor (v-1)/2 \rfloor = \lfloor (v-2)/2 \rfloor = \kappa - 1$. Therefore (3.17) and (3.20) gives

$$\sum_{l=1}^{\kappa} \hat{c}_{l,v+1} \lambda^l = \sum_{l=1}^{\kappa-1} v\hat{c}_{l,v-1} \lambda^{l+1} + \sum_{l=1}^{\kappa} l\hat{c}_{l,v} \lambda^l = \sum_{l=2}^{\kappa} (v\hat{c}_{l-1,v-1} + l\hat{c}_{l,v}) \lambda^l + \lambda.$$

For the even $v$ the property (3.19) follows. Let now $v$ is odd, then $\lfloor (v-1)/2 \rfloor = \kappa$ and $\lfloor (v+1)/2 \rfloor = \kappa + 1$. From (3.17) and (3.20) obtain

$$\sum_{l=1}^{\kappa+1} \hat{c}_{l,v+1} \lambda^l = \sum_{l=1}^{\kappa} v\hat{c}_{l,v-1} \lambda^{l+1} + \sum_{l=1}^{\kappa+1} l\hat{c}_{l,v} \lambda^l = \sum_{l=2}^{\kappa+1} v\hat{c}_{l-1,v-1} \lambda^l + \sum_{l=1}^{\kappa} l\hat{c}_{l,v} \lambda^l$$

$$= \sum_{l=2}^{\kappa} (v\hat{c}_{l-1,v-1} + l\hat{c}_{l,v}) \lambda^l + v\hat{c}_{\kappa,v-1} \lambda^{\kappa+1} + \lambda.$$

This proves (3.19) for the odd $v$. Lemma 3.2 is proved completely.

Using (3.15) we obtain $\mu_v(\lambda) < \lambda \frac{d}{d\lambda} \mu_v(\lambda) \leq \lfloor v/2 \rfloor \mu_v(\lambda)$. This together with (3.20) imply

$$\mu_v(\lambda) < \mu_{v+1}(\lambda) \leq (v\lambda + \lfloor v/2 \rfloor) \mu_v(\lambda). \qquad (3.21)$$

Use formula (3.20) for $\mu_{v+2}(\lambda)$, next again apply (3.20) to the derivative of $\mu_{v+1}(\lambda)$, after this use the inequalities $\lambda \frac{d}{d\lambda} \mu_{v-1}(\lambda) \leq \mu_v(\lambda)$ (because $l\hat{c}_{l,v-1} \leq \hat{c}_{l,v}$ (see (3.19)), $\lambda \frac{d}{d\lambda} \mu_v(\lambda) \leq \lfloor v/2 \rfloor \mu_v(\lambda)$ and

$\lambda^2 \frac{d^2}{d\lambda^2} \mu_v(\lambda) \leq \lfloor v/2 \rfloor^2 \mu_v(\lambda)$. Then one can observe that for any integer $v \geq 0$

$$(v+1)\lambda \mu_v(\lambda) < \mu_{v+2}(\lambda) \leq 2v(\lambda + v)\mu_v(\lambda). \qquad (3.22)$$

In the sequel some notations regarding to the general statistic $R_N$ (see (3.1)) are useful:

$$A_N = \sum_{m=1}^{N} E f_m(\xi_m), \quad \gamma_n = n^{-1} \sum_{m=1}^{N} \text{cov}(f_m(\xi_m), \xi_m),$$

$$\tilde{\sigma}_N^2 = \sum_{m=1}^{N} \text{Var} f_m(\xi_m), \quad \sigma_N^2 = \tilde{\sigma}_N^2 - n\gamma_n^2. \tag{3.23}$$

Note that under the conditions (3.2) and (3.3)

$$ER_N = A_N + o(N) \text{ and } \text{Var} R_N = \sigma_N^2(1+o(1)). \tag{3.24}$$

We notice that formula for the $\sigma_N^2$ and $\tilde{\sigma}_N^2$ given in (2.1) are the special case of (3.23) corresponding to the chi-square statistic; that is why the same symbols are used.

**Proof of Theorem 2.1**. The chi-square statistic is (3.1) kind with $f_m(x) = (x - np_m)^2 / np_m$, hence $A_N = N$, $\sigma_N^2$ and $\tilde{\sigma}_N^2$ as in (2.1). Set $\tilde{\xi}_m = (\xi_m - np_m)/\sqrt{np_m}$. We have $E\tilde{\xi}_m^{2s}(\xi_m - np_m)^2 \leq 4s(np_m + 2s)E\tilde{\xi}_m^{2s}$, since (3.22), hence the conditions (3.2) and (3.3) are fulfilled with $b=1$ and $a=1$ if $np_m \geq 2s$, and with $b=0$ and $a=2$ if $np_m < 2s$. It is evident that the condition (3.4) is fulfilled for the $k = o(p_{\max}^{-1})$. Therefore, from Main Lemma and Remark 3.1 it follows that

$$\mathcal{C}_k(\chi_N^2) = \mathcal{C}_k\left(\sum_{m=1}^{N} \tilde{\xi}_m^2\right)(1+o(1)), \tag{3.25}$$

for every integer $k$ such that $3 \leq k = o\left(\min\left(p_{\max}^{-1/2}, n^{1/3}\right)\right)$.

We would remind the notation $\nabla_n = \max\, 1, (np_{\min})^{-1}$. Since $\sqrt{2\pi m}\, m^m e^{-m} \leq m! \leq e\sqrt{m}\, m^m e^{-m}$ one can see that $(2k)! \leq 2^{2k} e/\pi \sqrt{2k}\, (k!)^2$. Due to this, putting $\varsigma_m = (\tilde{\xi}_m^2 - E\tilde{\xi}_m^2)/\sqrt{\text{Var}\tilde{\xi}_m^2}$ and using Lemma 3.2 we obtain $\left|E\varsigma_m^k\right| \leq 2^k (2k)!(\text{Var}\tilde{\xi}_m^2)^{-k/2} \sum_{l=0}^{k-1}(np_m)^{-l} \leq (k!)^2 \left(2^7 \nabla_n / \sqrt{\text{Var}\tilde{\xi}_m^2}\right)^{k-2}$ for $k \geq 3$, since

$\sum_{l=0}^{k-1}(np_m)^{-l} \leq k \nabla_n^{k-2}$ and $\text{Var}\tilde{\xi}_m^2 = (np_m)^{-1} + 2$. Hence, by Assertion 1

$\left|\mathcal{C}_k(\varsigma_m)\right| \leq (k!)^2 \left(2^8 \nabla_n / \sqrt{\text{Var}\tilde{\xi}_m^2}\right)^{k-2}$. That is $\left|\mathcal{C}_k(\tilde{\xi}_m^2)\right| \leq (k!)^2 \left(2^8 \nabla_n\right)^{k-2} \text{Var}\tilde{\xi}_m^2$, for $k \geq 3$. Therefore

$$\left|\mathcal{C}_k\left(\sum_{m=1}^{N}\tilde{\xi}_m^2\right)\right| \leq (k!)^2 \left(2^8 \nabla_n\right)^{k-2} \tilde{\sigma}_N^2, \tag{3.26}$$

because $\tilde{\xi}_m$, $m=1,2,\ldots$ are independent r.v.s, and $\tilde{\sigma}_N^2 = \sum_{m=1}^{N} \text{Var}\tilde{\xi}_m^2$, see notation (2.1). Now (3.25), (3.26) and the fact that $\sigma_N^2 = \text{Var}\chi_N^2(1+o(1))$, see (3.23), imply

$$\left|\mathcal{C}_k\left(\chi_N^2 / \sqrt{\text{Var}\chi_N^2}\right)\right| \leq \mathcal{C}_k\left(\sum_{m=1}^{N}\tilde{\xi}_m^2 / \sqrt{\text{Var}\chi_N^2}\right) \leq (k!)^2 \left(\frac{2^{13} \nabla_n \tilde{\sigma}_N^2}{\sigma_N^3}\right)^{k-2},$$

for all $k$ satisfying $3 \leq k = o\left(\min\left(p_{\max}^{-1/2}, n^{1/3}\right)\right)$. Theorem 2.1 follows from this and Assertion 2 with $\Delta = 2^{13} \nabla_n \tilde{\sigma}_N^2 / \sigma_N^3$ and $\tilde{k} = \min\left(p_{\max}^{-1/2}, n^{1/3}\right)$.

**Proof of Corollary 2.1**. Under conditions of Corollary 2.1 we have $\tilde{\sigma}_N^2\left(1 - (2c_0)^{-1}\right) \leq \sigma_N^2 \leq \tilde{\sigma}_N^2$ and $\nabla_n \leq \max(1, 2/c_1)$, also $p_{\max}^{-1/4} \geq c_2^{-1/4} N^{1/6}$. Corollary 2.1 follows.

**Proof of Corollary 2.2** is straightforward since in this case $\nabla_n = \max(1, \lambda_n^{-1})$,
$$\sigma_N^2 = N\left(2 + d^2 \delta^2(n) \lambda_n^{-1}(1 + o(1))\right) \text{ and } \tilde{\sigma}_N^2 = 2N\left(1 + (2\lambda_n)^{-1}(1 + o(1))\right) \leq 2N \nabla_n.$$

**Proof of Theorem 2.2.** Observe that the $\Lambda_N$ is $R_N$ kind statistic with $f_m(x) = x \ln(x/np_m)$ (we assume that $0 \cdot \ln 0 = 0$). Write
$$f_m(\xi_m) = np_m\left(1 + \hat{\xi}_m\right)\ln\left(1 + \hat{\xi}_m\right), \tag{3.27}$$
where $\hat{\xi}_m = (\xi_m - np_m)/np_m > -1$. Put $g_m(x) = f_m(x) - \gamma_n(x - np_m)$. It is evident that $R_N = \sum_{m=1}^{N} g_m(\eta_m)$, hence it is sufficient to check the conditions of Main Lemma for $g_m(\xi_m)$. We shall use the inequality
$$0 \leq (1+x)\ln(1+x) - \left(x + \frac{1}{2}x^2 - \frac{1}{6}x^3 + \ldots - \frac{x^{2l+1}}{2l(2l+1)}\right) \leq \frac{x^{2l+2}}{(2l+1)}, \tag{3.28}$$
$x > -1$ and integer $l \geq 0$, to find, see denotes (3.23),
$$A_N = \frac{N}{2} + \frac{1}{12}\sum_{m=1}^{N}\frac{1}{np_m}\left(1 + O\left(\frac{1}{np_m}\right)\right) = \frac{N}{2}\left(1 + O\left(\frac{1}{\lambda_n}\right)\right), \tag{3.29}$$
$$\gamma_n = 1 - \frac{1}{12n}\sum_{m=1}^{N}\frac{1}{np_m}\left(1 + O\left(\frac{1}{np_m}\right)\right) = 1 + O\left(\lambda_n^{-2}\right), \tag{3.30}$$
$$\sigma_N^2 = \sum_{m=1}^{N} Var\, g_m(\xi_m) = \frac{N}{2} - \frac{1}{6}\sum_{m=1}^{N}\frac{1}{np_m}\left(1 + O\left(\frac{1}{np_m}\right)\right) = \frac{N}{2}\left(1 + O\left(\frac{1}{\lambda_n}\right)\right). \tag{3.31}$$

From (3.28) we have $x \leq (1+x)\ln(1+x) \leq x + x^2/2$, hence
$$-\frac{c}{\lambda_n} \leq np_m(1-\gamma_n)\hat{\xi}_m \leq g_m(\xi_m) \leq np_m\left((1-\gamma_n)\hat{\xi}_m + \hat{\xi}_m^2\right), \tag{3.32}$$
since $\hat{\xi}_m \geq -1$, $np_m \geq c_1 \lambda_n$ and $1 - \gamma_n \leq (c_0 \lambda_n^2)^{-1}$ by (3.30). From Lemma 3.2 we obtain $E\hat{\xi}_m^\ell = O\left(\ell!(np_m)^{-\ell/2}\right)$ for any integer $\ell \geq 1$ (to see this for an even $\ell \geq 1$ we can directly use Lemma 3.2, whereas for an odd $\ell \geq 1$ we first use the inequality $E|\hat{\xi}_m^\ell| \leq \left(E\hat{\xi}_m^{\ell+1}\right)^{\ell/(\ell+1)}$ and then apply Lemma 3.2). Therefore, it is not hard to understand that the conditions (3.2) and (3.3) with $a=1$ and $b=1$ are satisfied. Thus, Main Lemma gives: for all integer $k$ such that $3 \leq k = o\left(p_{\max}^{-1/2}\right)$

$$\mathcal{C}_k(\Lambda_N) = \mathcal{C}_k\left(\sum_{m=1}^{N} g_m(\xi_m)\right)(1+o(1)). \tag{3.33}$$

On the other hand using (3.32), inequality $(2k)! \leq 2^{2k} e/\pi \sqrt{2k} \, (k!)^2$ and the condition (2.4) we obtain

$$E|g_m(\xi_m)|^k \leq 2^k (np_m)^k \left((1-\gamma_n)^k E|\hat{\xi}_m|^k + E\hat{\xi}_m^{2k}\right) \leq 2^{k+1}\left((2k)! + \frac{k!(np_m)^{k/2}}{c_0 \lambda_n^{2k}}\right) \leq (k!)^2 \left(2^{11}\right)^{k-2}.$$

Therefore, Assertion 1 imply $\left|\mathcal{C}_k(g_m(\xi_m))\right| \leq (k!)^2 \left(2^{12}\right)^{k-2}$, because $\mathrm{Var}\, g_m(\xi_m) = 2^{-1}(1+o(1))$. Hence

$$\left|\mathcal{C}_k\left(\sigma_N^{-1}\sum_{m=1}^{N} g_m(\xi_m)\right)\right| \leq (k!)^2 \left(2^{11} N^{-1/2}\right)^{k-2},$$

since (3.29) and (3.31). Thus from (3.33) taking into account that $\sigma_N^2 = \sum_{m=1}^{N} \mathrm{Var}\, g_m(\xi_m)$ we obtain

$$\left|\mathcal{C}_k(\sigma_N^{-1}\Lambda_N)\right| \leq (k!)^2 \left(2^{12} N^{-1/2}\right)^{k-2},$$

for all integer $k$ such that $3 \leq k = o\left(p_{\max}^{-1/2}\right)$. Noting that $\mathrm{Var}\,\Lambda_N = \sigma_N^2(1+o(1))$, see (3.24), Theorem 2.2 follows from Assertion 2 with $\tilde{k} = o\left(p_{\max}^{-1/2}\right)$ and $\Delta = 2^{12} N^{-1/2}$.

**Proof of Corollary 2.3** is straightforward.

**Appendix.** Below $\xi$ is a r.v. with $E\xi = 0$, $\mathrm{Var}\,\xi = 1$, $\mathcal{C}_k(\xi)$ and $\alpha_k(\xi)$ cumulant and moment of $k$ th order respectively of the r.v. $\xi$. The following conditions and assertions are taken from the book by Saulis and Statulevicius (1991).

**Bernstein condition** $(B_\nu)$: there exists the constants $\nu \geq 0$ and $B > 0$ such that

$$|\alpha_k(\xi)| \leq (k!)^{\nu+1} B^{k-2}, \text{ for all } k = 3, 4, \ldots.$$

**Statulevicius condition** $(S_\nu)$: there exists the constants $\nu \geq 0$ and $\Delta > 0$ such that

$$|\mathcal{C}_k(\xi)| \leq (k!)^{\nu+1} \Delta^{-(k-2)}, \text{ for all } k = 3, 4, \ldots.$$

**Assertion 1.** If $\xi$ satisfy condition $(B_\nu)$ then it also satisfy condition $(S_\nu)$ with $\Delta = (2B)^{-1}$.

**Assertion 2.** Let a r.v. $\xi$ depending on a parameter $\Delta$ satisfy condition $(S_1)$ (i.e. $\nu = 1$) for all $k$ such that $3 \leq k \leq \tilde{k}$. Then for all $x$ such that $0 \leq x \leq 12^{-1} \min\left(\sqrt{\tilde{k}}, \Delta^{1/3}\right)$ there exist a constant $c$ such that $P(\xi > x) = (1-\Phi(x))\left(1+C(x+1)\Delta^{-1/3}\right)$ and $P(\xi < -x) = \Phi(-x)\left(1+C(x+1)\Delta^{-1/3}\right)$.

In the book of Saulis and Statulevicius (1991) no assertion in the form of Assertion 2 is given. However, Assertion 2 can be proved by reasoning alike to those presented by Saulis and Statulevicius (1991) at page 35, in course of proof of Lemma 2.3, see relations (2.62), (2.65) and (2.66) and Lemma 2.2 of the book.

**Assertion 3**. One has

$$C_k(\xi) = k! \sum (-1)^{m_1+m_2+...+m_k-1} (m_1+m_2+...+m_k-1)! \prod_{l=1}^{k} \frac{1}{m_l!} \left( \frac{\alpha_l(\xi)}{l!} \right)^{m_l}$$

here $\sum$ is summation over all non-negative integer $m_1, m_2, ..., m_k$ such that $m_1 + 2m_2 + ... + km_k = k$.

**References.**